\documentclass{article}

\usepackage{amsmath,amsfonts,latexsym,theorem, 
}
\nonstopmode

\numberwithin{equation}{section}

\newtheorem{Theorem}{Theorem}[section]
\newtheorem{Definition}{Definition}[section]
\newtheorem{Proposition}{Proposition}[section]
\newtheorem{Lemma}{Lemma}[section]

\newenvironment{Proofc}[1]{\smallskip\par\noindent\textsc{#1}\quad}%
  {\hfill$\Box$\bigskip\par}
\newenvironment{Proof}{\begin{Proofc}{Proof}}{\end{Proofc}}

\newtheorem{Remark}{Remark}[section]

\def\b{\beta}

\def\G{\Gamma}

\def\e{\varepsilon}

\def\pd{\partial}
\def\half{\frac{1}{2}}

\newcommand{\cL}{{\mathcal L}}

\newcommand{\R}{{\mathbb R}}
\newcommand{\N}{{\mathbb N}}

\def\pd{\partial}

\begin{document}
\title{A note on  Kazdan-Warner equation on   networks}

\author{Fabio Camilli\footnotemark[1] \and  Claudio Marchi\footnotemark[2]}

\date{version: \today}
\maketitle

\footnotetext[1]{Dip. di Scienze di Base e Applicate per l'Ingegneria,   Sapienza  Universit{\`a}  di Roma, via Scarpa 16, 00161 Roma, Italy, ({\tt e-mail: fabio.camilli@uniroma1.it})}
\footnotetext[2]{Dip. di Ingegneria dell'Informazione, Universit\`a di Padova, via Gradenigo 6/B, 35131 Padova, Italy ({\tt claudio.marchi@unipd.it}).}
\begin{abstract}
We investigate the Kazdan-Warner  equation on a network. In this case,  the differential equation   is defined on each edge, while appropriate transition conditions of Kirchhoff type are prescribed  at the vertices. We show that the Kazdan-Warner theory extends to the present setting and we study also the critical case.
\end{abstract}
 \begin{description}
	\item [\textbf{MSC 2010}:] 35A15, 35J60,  35R02.
	\item [\textbf{ Keywords}:] Kazdan-Warner equation, network,  Kirchhoff condition.
\end{description}
%
%

\section{Introduction}\label{intro}
The Kazdan-Warner equation
\begin{equation}\label{intro_1}
\Delta u = c-he^u,	
\end{equation}
where $c$ is a constant and $h$ a given function, was introduced in \cite{kw} in connection with the problem of prescribing the Gaussian curvature of a compact manifold $M$. The solvability of \eqref{intro_1} depends on the sign of $c$. Let $\bar h$ denote the average of $h$ on $M$. In \cite{kw}, it  is shown   that
\begin{enumerate}
	\item if $c=0$ and $h\not\equiv 0$, then \eqref{intro_1} is solvable if and and only if $h$ changes sign and
	$\bar h<0$;
	\item if $c>0$,  then \eqref{intro_1} is solvable if and only if the set $\{h>0\}$ is not empty;
	\item if $c<0$,  if    \eqref{intro_1}  is solvable, then $\bar h<0$. For $\bar h<0$, there exists a constant $-\infty\le c(h)<0$ such that  \eqref{intro_1} is solvable for any $c\in (c(h),0)$ and not solvable for any $c<c(h)$. Moreover $c(h)=-\infty$ if and only if $h\leq 0$ in~$M$.
\end{enumerate}
If $c<0$,   $c=c(h)$ is not included in the previous cases and deserves a particular attention. It has been shown in \cite{cl} that, if $c(h)>-\infty$, then   \eqref{intro_1} can be also solved for $c=c(h)$.\par
The previous theory has been recently extended in \cite{ge,gly} to the case of a connected, finite graph.
Here the Laplacian is replaced by a finite difference operator, the so-called graph Laplacian, and most of the effort is to reproduce in a finite dimensional setting some crucial properties as the Maximum Principle and the Moser-Trudinger inequality.\par
An intermediate situation between a compact manifold and a finite graph is given by a network~$\G$, which is  given by a finite collection of vertices   connected by continuous non-self-intersecting edges. The differential equation \eqref{intro_1} is defined on each edge, while appropriate transition conditions of Kirchhoff type are prescribed  at the vertices. In this paper, we obtain the same conclusions of the manifold and finite graph cases,  showing that the Kazdan-Warner theory remains unchanged for different classes of manifolds, also non regular such as in the case of networks. To prove these results, we shall adapt the method by Kazdan Warner~\cite[Thm5.3]{kw} (see also~\cite[Thm2]{gly}) and, for the critical case, some techniques of~\cite{cl,ge} with some specific arguments for networks.\par
The paper is organized as follows. In Section \ref{sec;1}, we introduce some notations and preliminary results. In Sections \ref{sec;2}, \ref{sec;3} and \ref{sec;4}, we study respectively the cases $c=0$, $c>0$ and $c<0$. In Section \ref{sec;4}, we also discuss the critical case $c=c(h)$.
\section{Notations, definitions  and preliminary results}\label{sec;1}
A network~$\G=(V,E)$ is a finite collection of points $V:=\{v_i\}_{i\in I}$ in $\R^n$ connected by continuous non-self-intersecting edges $E:=\{e_j\}_{j\in J}$, where any two edges can only have intersection at a vertex.  For $i\in I$ we set \[Inc_i:=\{j\in J:\,e_j \,\text{is incident to}\,v_i\}.\]
A coordinate $\pi_j:[0,l_j]\to\R^n$, with $l_j>0$,  is chosen to parameterize $e_j$,   i.e.
$e_j:=\pi_j((0,l_j))$.  We  assume that $\G$ is finite and connected and we  denote with  $|\Gamma|$ the sum of lengths of the edges $e_j$, $j\in J$.
 \par
For a function $u:  \G\to\R$  we denote by $u_j:[0,l_j]\to \R$ the restriction of $u$ to $e_j$, i.e. $u(x)=u_j(y)$ for $x\in e_j$, $y=\pi_j^{-1}(x)\in (0,l_j)$. Given $v_i\in V$, we denote by $\pd_j u(v_i)$  the  oriented derivative at $v_i$ along the arc $e_j$ defined   by
\[\pd_j u (v_i)=\lim_{x\in e_j,\,x\to v_i}\frac{u(\pi_j^{-1}(x))-u(\pi_j^{-1}(v_i))}{|\pi^{-1}_j(x)-\pi^{-1}_j(v_i)|},\]
if the limit exists, where $\pi_j$ is the  parametrization of arc $e_j$. For a function $\phi:\G\to\R$ and $A\subset\G$, we  set
\[\int_A \phi(x)dx:=\sum_{j} \int_{(0,l_j)\cap\pi_j^{-1}(A)}\phi(r)dr.\]
A function  $u$  is said continuous on $\G$ if it is continuous  with respect to the subspace topology of $\G$, i.e. $u_j\in C([0,l_j])$ for any $j\in J$ and
$u_j(\pi_j^{-1}(v_i))=u_k(\pi_k^{-1}(v_i))$  for any $i\in I$, $j,k\in Inc_i$.\\
 We introduce some  functional spaces for functions defined on the network.
 The space $L^p(\G)$, $p\ge 1$,  consists of the functions that are measurable and $p$-integrable on each edge $e_j$, $j\in J$. We set
 \[\|f\|_{L^p  }:=\left(\sum_{j\in J}\|f_j\|^p_{L^p(e_j)}\right)^{1/p}.\]
The space $L^\infty(\G)$ consists of the functions that are measurable and bounded on each edge $e_j$, $j\in J$. We set
\[\|f\|_{L^\infty  }:=\sup_{j\in J}\|f_j\|_{L^\infty(e_j)}.\]
 The Sobolev space $W^{k,p}(\G)$, $k\in\N$ and $p\ge 1$,  consists of all continuous functions on $\G$ that belong  to $W^{k,p}(e_j)$ for each  $j\in J$. We set
\[\|f\|_{W  ^{k,p}}:=\left(\sum_{l=0}^k\|\pd^l f\|^p_{L^p  }\right)^{1/p}.\]
As usual we set $H^k(\G):=W^{k,2}(\G)$, $k\in\N$.
The space  $C^{k}(\G)$ for $k\in \N$ consists of all continuous functions on $\G$  that belongs to $ C^{k}(e_j)$ for $j\in J$. The space  $C^k(\G)$ is a Banach space  with the norm
\[\|f\|_{C^k}=\max_{\b\le k}\|\pd^\beta f\|_{L^\infty  }.\]
		

The following proposition gives a   Poincar\'e inequality for the network
\begin{Lemma}\label{lemma:poinc}
For every function $f\in H^1(\Gamma)$ with $\int_\Gamma f(x)dx=0$, there holds
\begin{itemize}
\item[(i)] $|f(x)|\leq \sqrt{|\Gamma|}\,\|\partial f\|_{L^2}$;
\item[(ii)] $\int_\Gamma f^2(x) dx\leq |\Gamma|^2 \int_\Gamma |\partial f(x)|^2 dx.$
\end{itemize}
\end{Lemma}
\begin{Proof}
By definition of $H^1$, the function $f$ is continuous on $\Gamma$, hence there exists a point $x_0\in\Gamma$ such that $f(x_0)=0$.
Since $\G$ is connected, for any point $x\in\Gamma$ there exists a path $\gamma:(0,r)\rightarrow \Gamma$ on the network such that $\gamma(0)=x_0$, $\gamma(r)=x$, $|\gamma'(s)|=1$ and $r\leq |\Gamma|$. Hence, we have
\[
|f(x)|=|f(x_0)+\int_0^r(f\circ\gamma)'(s)ds|\leq\int_0^r|\partial f(\gamma(s))|ds\leq \sqrt{r}\|\partial f\|_{L^2(\gamma)}\leq \sqrt{|\Gamma|}\|\partial f\|_{L^2}.
\]
We deduce that
\[
\int_\Gamma f^2(x)dx\leq  \int_\Gamma |\Gamma|\|\partial f\|_{L^2}^2dx\leq |\Gamma|^2 \|\partial f\|_{L^2}^2.
\]
\end{Proof}
We also give an analogous of the Trudinger-Moser inequality for the networks. 
\begin{Lemma}\label{lemma:moser}
For any $\beta,\delta\in\R$ with~$\delta>0$, there exists a constant $C$ (depending only on $\beta$, $\delta$ and the network) such that, for all functions  $f\in H^1(\Gamma)$ with $\int_\Gamma |\partial f|^2\leq \delta$ and $\int_\Gamma f =0$, there holds
\[
\int_\Gamma e^{\beta f^2(x)}\,dx \leq C.
\]
\end{Lemma}
\begin{Proof}
We adapt the arguments of ~\cite[Lemma7]{gly}. The case $\beta \leq 0$ is obvious because $\Gamma$ has a bounded total length. Fix $\beta>0$ and consider a function~$f$ as in the statement. By Lemma~\eqref{lemma:poinc}-(i) and by the assumption $\|\partial f\|_2^2\leq \delta$, we have
\[
\int_\Gamma e^{\beta f^2(x)}\,dx \leq \int_\Gamma e^{\beta|\Gamma| \|\partial f\|_2^2}\,dx \leq e^{\beta |\Gamma| \delta}|\Gamma|.
\]

\end{Proof}

We consider the Kazdan-Warner equation on the network $\G$
\begin{equation}\label{KW}
\left\{\begin{array}{ll}
\partial^2 u = c-he^u,&\qquad \textrm{if } x\in e_j,\, j\in J,\\[4pt]
u_j(v_i)=u_k(v_i),&\qquad j,k\in Inc_i,\, v_i\in V,\\[4pt]
\sum_{j\in Inc_i} \partial_j u(v_i)=0,&\qquad v_i\in V,
\end{array}\right.
\end{equation}
where $c$ is given constant and $h$ is a continuous function on~$\Gamma$. Note that the Kazdan-Warner equation is defined on each edge, while at the vertices we impose the continuity of $u$ and the Kirchhoff condition, a classical condition for differential equations defined on networks (see \cite{mu,pb}).
\begin{Definition}
\begin{itemize}
\item[(a)] A strong solution to problem~\eqref{KW} is a function $u\in C^2(\Gamma)$ which satisfies~\eqref{KW} in a pointwise manner.
\item[(b)] A weak solution to problem~\eqref{KW} is a function $u\in H^1(\Gamma)$ such that
\begin{equation}\label{KW_weak}
\int_\Gamma\partial u \partial \phi\, dx= -c \int_\Gamma \phi\, dx +\int_\Gamma he^u \phi\, dx \qquad \forall \phi\in H^1(\Gamma).
\end{equation}
\end{itemize}
\end{Definition}
\begin{Remark}\label{rmk:weak-strong}
One can easily check that, if $u\in C^2(\Gamma)$ is a weak solution of~\eqref{KW}, then it is also a strong solution.
Moreover,  any weak solution of~\eqref{KW} is also a strong solution. Actually, a weak solution~$u$ fulfills $\partial^2 u=c-he^u$ in distributional sense inside each edge~$e_j$. The right hand side of this equality is continuous, hence, by standard theory, $u\in C^2(e_j)$ for every~$j\in J$. Being a weak solution, $u$ also belongs to~$H^1(\Gamma)$; in conclusion $u\in C^2(\Gamma)$.  
\end{Remark}
In the next three sections, we discuss the solvability of \eqref{KW} in the cases $c=0$, $c>0$ and $c<0$.
\section{The Kazdan-Warner equation with case $c=0$}
\label{sec;2}
\begin{Theorem}\label{thm:=0}
	Assume $c=0$ and $h\not\equiv 0$. Then problem~\eqref{KW} has a solution~$u$ if and only if~$h$ changes sign and $\int_\Gamma h<0$.
\end{Theorem}
\begin{Proof}
Assume that $u$ is a solution to problem~\eqref{KW} with $c=0$. We note that the hypothesis $h\not\equiv 0$ prevents~$u$ to be constant. 
We multiply the differential equation in~\eqref{KW} by $\phi\equiv 1$ and integrate on~$\Gamma$; taking advantage of the Kirchhoff condition, we get $\int_\Gamma he^udx=0$ which implies that $h$ must change sign.
Multiplying $e^{-u}\partial^2 u = -h$  by $\phi\equiv 1$ and integrating on~$\Gamma$, we get
\[
\int_\Gamma(\partial u)^2e^{-u}\, dx +\sum_{i\in I}\sum_{j\in Inc_i}e^{-u(v_i)}\partial_j u(v_i) = -\int_\Gamma h\, dx.
\]
Taking advantage of the Kirchhoff condition and of the continuity of $u$ at each vertex, we obtain
\[
\int_\Gamma(\partial u)^2e^{-u}\, dx = -\int_\Gamma h\, dx.
\]
Since $u$ cannot be constant, we deduce $\int_\Gamma hdx<0$.\par
Conversely, we   prove that, for any $h$ which changes sign and satisfies $\int_\Gamma h<0$,  there exists a solution to~\eqref{KW}. We define the set
\[
B:=\left\{v\in H^1(\Gamma)\mid \int_\Gamma he^vdx=0,\, \int_\Gamma vdx=0\right\}.
\]
We claim that $B$ is not empty. Since $h$ changes sign, there exists a point $x_0\in\Gamma$ such that $h(x_0)>0$. By the continuity of~$h$, without any loss of generality, we can assume $x_0\in e_{\bar \j}$ for some $\bar \j\in J$; namely, there exist $\bar\j\in J$ and $y_0\in (0,l_{\bar \j})$ such that $h_{\bar \j}(y_0)>0$.
Moreover, still by the continuity of~$h$, there exists~$\e>0$ such that $(y_0-\e,y_0+\e)\subset(0,l_{\bar \j})$ and  $h_{\bar \j}(y)>h_{\bar \j}(y_0)/2$ for all $y\in (y_0-\e,y_0+\e)$. Consider a function~$w\in C^2(\Gamma)$ such that: $w_{\bar \j}(y)=1$ if $y\in (y_0-\e/2,y_0+\e/2)$, $w_{\bar \j}(y)=0$ if $y\notin (y_0-\e,y_0+\e)$ and $w_{j}\equiv 0$ if $j\in J\setminus\{\bar\j\}$. For $\ell>0$, the function $w_\ell(\cdot):=\ell w (\cdot)$ fulfills
\begin{eqnarray}\notag
  \int_\Gamma he^{w_\ell}\, dx &=& \int_{e_{\bar \j}}h e^{w_\ell}\, dx +\sum_{j\in J\setminus\{\bar\j\}}\int_{e_j} h e^{w_\ell}\, dx\geq \int_{y_0-\e/2}^{y_0+\e/2}h(y)e^{w_\ell(y)}\, dy -\int_\Gamma |h| \, dx\\ \label{ell}
&\geq & 
\frac{\e h_{\bar \j}(y_0)e^\ell}{2} -\int_\Gamma |h| \, dx>0
\end{eqnarray}
provided that $\ell$ is sufficiently large. On the other hand, for $\ell=0$ we have $w_0(x)\equiv 0$ and, by  assumptions,
\begin{eqnarray*}
\int_\Gamma he^{w_0(x)}\, dx = \int_\Gamma h \, dx<0.
\end{eqnarray*}
Therefore there exists $\ell_0>0$ such that $\int_\Gamma he^{w_{\ell_0}}=0$. Hence
the function $\hat w(\cdot):=w_{\ell_0}(\cdot)- \int_\Gamma w_{\ell_0}/|\Gamma|$ belongs to $B$ and the claim is proved.\par
Consider the functional
\[
{\mathcal J}(v):=\frac12\int_\Gamma |\partial v|^2\, dx, \qquad \forall v\in B.
\]
Let $\{v_n\}_{n\in\N}$ be a minimizing sequence for~${\mathcal J}$, i.e. $\lim_{n\to+\infty} {\mathcal J}(v_n)=\inf_B {\mathcal J}$. By Lemma~\ref{lemma:poinc}-(ii), possibly passing to a subsequence, we have that the functions $v_n$ are uniformly bounded in~$H^1(\Gamma)$. We deduce that there exists $\bar u\in H^1(\Gamma)$ such that, as $n\to+\infty$, $v_n\rightharpoonup \bar u$ weakly in~$H^1(\Gamma)$ and $v_n\to \bar u$ uniformly on ~$\Gamma$. In particular, we get that $\bar u$ belongs to~$B$ and it is a minimizer of~${\mathcal J}$ on~$B$.\par
We claim that~$\bar u$ is a strong solution to problem~\eqref{KW}. Actually, by standard Lagrangian multiplier theory, there exist $\lambda,\mu\in\R$ such that
\begin{eqnarray*}
0&=&\frac{d}{dt}\left.\left(
{\mathcal J}(\bar u+t\phi)-\lambda \int_\Gamma he^{\bar u+t\phi}\, dx -\mu \int_\Gamma (\bar u+t\phi) \, dx \right)\right|_{t=0}\\
&=& \int_\Gamma \partial \bar u\partial \phi\, dx-\lambda \int_\Gamma he^{\bar u}\phi\, dx -\mu \int_\Gamma \phi \, dx,
\end{eqnarray*}
for every $\phi\in H^1(\Gamma)$.
Choosing $\phi\equiv 1$, since $\bar u\in B$, we get $\mu=0$. Arguing as in Remark~\ref{rmk:weak-strong}, inside each edge~$e_j$, there holds $\partial ^2\bar u+\lambda h e^{\bar u}=0$ in distributional sense. By the continuity of~$\bar u$, we infer that~$\bar u\in C^2(e_j)$ and, since $\bar u\in H^1(\Gamma)$, also that $\bar u\in C^2(\Gamma)$. Moreover, $\bar u$ is a strong solution to 
\begin{equation}\label{zero_1}
\left\{\begin{array}{ll}
\partial^2 \bar u = -\lambda he^{\bar u}&\qquad \textrm{if } x\in e_j,\, j\in J\\
\sum_{j\in Inc_i} \partial_j \bar u(v_i)=0&\qquad v_i\in V\\
\bar u_j(v_i)=\bar u_k(v_i)&\qquad j,k\in Inc_i,\, v_i\in V.
\end{array}\right.
\end{equation}
We claim $\lambda >0$. 
The function~$\bar u$ also solves $e^{-\bar u}\partial^2 \bar u = -\lambda h$; integrating this relation, by Kirchhoff and continuity conditions, we get
\[
\int_\Gamma(\partial \bar u)^2e^{-\bar u}\, dx = -\lambda\int_\Gamma h\, dx.
\]
Let us first prove that the left hand side of this equality is positive. We proceed by contradiction assuming $\int_\Gamma(\partial \bar u)^2e^{-\bar u}=0$. Hence, $\partial \bar u\equiv 0$ and, in particular, $\bar u$ is constant. Since $\bar u\in B$, we get $e^{\bar u}\int_\Gamma h=0$ contradicting the assumption $\int_\Gamma h<0$.
Therefore, the left hand side in the last equality is positive; again by virtue of $\int_\Gamma h<0$, the constant~$\lambda$ must be positive.
Finally, the function $u(\cdot):=\bar u(\cdot)+c$ with $c:=\log(\lambda)$ is a strong solution to~\eqref{KW}.
\end{Proof}    
\section{The Kazdan-Warner equation with case $c>0$}
\label{sec;3}
\begin{Theorem}\label{thm:>0}
Assume $c>0$. Then problem~\eqref{KW} has a solution~$u$ if and only if~$h$ is positive somewhere.
\end{Theorem}
\begin{Proof}
Assume that $u$ is  a solution of \eqref{KW}; choosing $\phi\equiv 1$ as test function in \eqref{KW_weak}, we get $\int_\Gamma h e^u = c|\Gamma|>0$. Hence, $\{x\in \Gamma\mid h(x)>0\}\neq \emptyset$.

Conversely, for any~$h\in C^0(\Gamma)$ with $\{h>0\}\neq \emptyset$, we prove that the problem~\eqref{KW} admits at least one solution. To this end, it is expedient to introduce the set
\[
B:=\left\{
v\in H^1(\Gamma)\mid \int_\Gamma h e^u\, dx = c|\Gamma|\right\}.
\]
We claim that~$B$ is not empty.
For~$\ell\geq 0$, we introduce the function~$w_\ell$ as in the proof of Theorem~\ref{thm:=0} while, for~$\ell\leq 0$, we set~$\bar w_\ell\equiv \ell$. Since~$w_0\equiv \bar w_0$, the function
\[
g(\ell):=\left\{\begin{array}{ll}
\int_\Gamma h e^{w_\ell}\, dx &\textrm{ if }\ell\geq 0\\
\int_\Gamma h e^{\bar w_\ell}\, dx &\textrm{ if }\ell< 0
\end{array}\right.
\]
is well defined and continuous, it fulfills $\lim_{\ell \to +\infty}g(\ell)=+\infty$ (by virtue of the estimate~\eqref{ell}) and $\lim_{\ell \to -\infty}g(\ell)=\lim_{\ell \to -\infty} e^{\ell}\int_\Gamma h=0$. Hence, there exists~$\bar \ell\in \R$ such that $g(\bar \ell)=c|\Gamma|$ namely,~$B\neq \emptyset$.

We consider the functional
\[
{\mathcal J}(u):=\frac12 \int_\Gamma |\partial u|^2\, dx + c \int_\Gamma u \, dx, \qquad \forall u\in B.
\]
As a first step, let us prove that ${\mathcal J}$ is bounded from below. To this end, for any~$u\in B$, we set~$\bar u:=\int_\Gamma u/|\Gamma|$ and $v:=u-\bar u$. Note $\int_\Gamma v=0$ and $\partial v\equiv \partial u$. Since~$u\in B$, it holds $\int_\Gamma h e^{v}\, dx = c|\Gamma|e^{-\bar u}$ which implies $\bar u=\log(c|\Gamma|)-\log \left(\int_\Gamma h e^{v}\, dx\right)$; replacing this equality in the definition of~${\mathcal J}$, we get
\begin{equation}\label{c>0;3}
{\mathcal J}(u)= \frac12 \|\partial u\|_2^2 + c|\Gamma|\log(c|\Gamma|)-c|\Gamma| \log \left(\int_\Gamma h e^{v}\, dx\right).
\end{equation}
Let us now estimate $\int_\Gamma h e^{v}$; if $v$ is constant then, by $\int_\Gamma v=0$, it must be $v\equiv 0$ and, in particular $\int_\Gamma h e^{v}=\int_\Gamma h$. For~$v$ nonconstant, it is expedient to introduce the function~$\tilde v:=v/\| \partial v\|_2$ which verifies: $\tilde v\in H^1(\Gamma)$, $\int_\Gamma \tilde v=0$ and~$\|\partial \tilde v\|_2=1$.
Lemma~\ref{lemma:poinc}-(ii) and Lemma~\ref{lemma:moser} guarantee that, for any~$\beta\in\R$, there exists a constant~$K_\beta$ (depending only on~$\beta$) such that
\begin{equation*}
\|\tilde v\|_2\leq |\Gamma|, \qquad
\int_\Gamma e^{\beta \tilde v^2(x)}\, dx\leq K_\beta.
\end{equation*}
For every~$\e$ positive, for $\beta_\e:=1/(4\e)$, there holds
\[
\int_\Gamma h e^{v}\, dx\leq \|h\|_\infty \int_\Gamma e^{\e \|\partial v\|_2^2 +\frac{v^2}{4\e \|\partial v\|_2^2}}\, dx\leq \|h\|_\infty e^{\e \|\partial v\|_2^2 } K_{\beta_\e}.
\]
Replacing this estimate in~\eqref{c>0;3}, we obtain
\begin{equation*}
{\mathcal J}(u)\geq \frac12 \|\partial u\|_2^2  + c|\Gamma|\left[\log(c|\Gamma|) - \e \|\partial u\|_2^2 -\log \left(\|h\|_\infty K_{\beta_\e}\right)\right]
\end{equation*}
and, in particular, for~$\e_0:=\frac{1}{4 c|\Gamma|}$,
\begin{equation}\label{c>0;4}
{\mathcal J}(u)\geq \frac14 \|\partial u\|_2^2  + c|\Gamma|\left[\log(c|\Gamma|) -\log \left(\|h\|_\infty K_{\beta_{\e_0}}\right)\right].
\end{equation}
Hence, the proof that~${\mathcal J}$ is bounded from below is accomplished.

Let $\{u_n\}_{n\in\N}$ be a minimizing sequence for~${\mathcal J}$; set $\bar u_n:= \int_\Gamma u_n/|\Gamma|$ and $v_n:=u_n-\bar u_n$; hence $\partial u_n\equiv \partial v_n$ and, by estimate~\eqref{c>0;4},  $\partial v_n$ is bounded in $L^2(\Gamma)$, uniformly in~$n$. By Lemma~\eqref{lemma:poinc}-(ii), also $v_n$ is uniformly bounded in $L^2(\Gamma)$ and, therefore, the functions~$v_n$ are uniformly bounded in~$H^1(\Gamma)$.
Moreover, by the definition of~${\mathcal J}$, we get that $\int_\Gamma u_n$ are uniformly bounded and consequently also $\bar u_n$ are uniformly bounded. Being $u_n=v_n+\bar u_n$, also the functions~$u_n$ are uniformly bounded in~$H^1(\Gamma)$. Possibly passing to a subsequence, there exists~$u\in H^1(\Gamma)$ such that, as $n\to+\infty$,  $u_n\rightharpoonup u$ in the weak topology of~$H^1(\Gamma)$, $u_n\to u$ uniformly, $u\in B$ and ${\mathcal J}(u)=\min_B {\mathcal J}$.

We claim that $u$ is a solution to~\eqref{KW}. By standard Lagrangian theory, there exists~$\lambda\in \R$ such that, for every $\phi\in H^1(\Gamma)$,
\begin{eqnarray}\notag
0&=&\frac{d}{dt}\left.\left(
\int_\Gamma\frac{\partial (u+t\phi)^2}{2}\, dx+c\int_\Gamma (u+t\phi)\, dx -\lambda \left(c|\Gamma|-\int_\Gamma he^{u+t\phi}\, dx\right) \right)\right|_{t=0}\\ \label{c>0;5}
&=& \int_\Gamma \partial u\partial \phi\, dx+c\int_\Gamma \phi\, dx -\lambda \int_\Gamma he^{u}\phi\, dx.
\end{eqnarray}
Choosing $\phi\equiv 1$, we get $c|\Gamma|=\lambda \int_\Gamma he^{u}$; since~$u\in B$, we get~$\lambda =1$. In conclusion, relation~\eqref{c>0;5} with~$\lambda =1$ is equivalent to the definition of weak solution to~\eqref{KW}.
\end{Proof}
\section{The Kazdan-Warner equation with case $c<0$}
\label{sec;4}
\begin{Theorem}\label{thm:<0}
Assume $c<0$. Then
\begin{itemize}
\item[(i)] If \eqref{KW} has a solution, then $\int_\Gamma h<0$.
\item[(ii)] If $\int_\Gamma h<0$, then there exists a constant $c(h)\in [-\infty,0)$ such that \eqref{KW} has a solution for any $c(h)<c<0$  and  no solution for $c<c(h)$.
\item[(iii)] For $\int_\Gamma h<0$, let $c(h)$ be defined as in {\it (ii)}. Then, $c(h)=-\infty$ if and only if $h\le 0$ in $\G$.

\end{itemize}	
\end{Theorem}
We introduce the definition of upper and lower solution to \eqref{KW}.
\begin{Definition}
A function $u \in C^2(\Gamma)$ is said to be
a lower (respectively, an upper) solution of \eqref{KW} if
\[
\left\{\begin{array}{ll}
\partial^2 u  - c+he^u \ge 0& \textrm{if } x\in e_j,\, j\in J,\\
\sum_{j\in Inc_i} \partial_j u (v_i)\ge 0& v_i\in V,
\end{array}\right.
\ \left(\textrm{resp.,}
\left\{\begin{array}{ll}
\partial^2 u  - c+he^u \le 0& \textrm{if } x\in e_j,\, j\in J,\\
\sum_{j\in Inc_i} \partial_j u (v_i)\le 0& v_i\in V
\end{array}
\right.\right).\]
\end{Definition}
In order to prove Theorem \ref{thm:<0}, we need some preliminary results.
\begin{Lemma}\label{c<0;L1}
If there exist a lower solution $u_-$ and an upper solution $u_+$ of \eqref{KW}	such that $u_-\le u_+$, then there
there exists a solution $u$ of \eqref{KW} such that $u_-\le u\le u_+$.
\end{Lemma}
\begin{Proof}
Set $k_1(x)=\max\{1, -h(x)\}$ and $k(x)=k_1(x)e^{u_+(x)}$ and consider the sequence of function $\{u_n\}_{n\in\N}$ defined inductively as $u_0=u_+$ and $u_n$   the solution of 
\begin{equation}\label{c<0;1}
\left\{\begin{array}{ll}
\cL u_{n+1}=f(x,u_n)-ku_n&\qquad \textrm{if } x\in e_j,\, j\in J\\
\sum_{j\in Inc_i} \partial_j u (v_i)=0&\qquad v_i\in V
\end{array}\right.
\end{equation}
where $\cL u=\partial^2 u -k u$ and $f(x,u)=c-h(x)e^u$. We first observe that the sequence $\{u_n\}_{n\in\N}$ is well defined: indeed, since $k(x)\ge e^{-\|u_+\|_\infty}$,    \eqref{c<0;1} admits a unique strong  solution   $u_n$ for any $n\in \N$ (see \cite[Prop.10]{cm}). Moreover, we claim that  
\begin{equation}\label{c<0;2}
u_-\le u_{n+1}\le u_n\le u_+, \quad \text{for any $n\in \N$.}
\end{equation}
Since 
\[
\left\{\begin{array}{ll}
\cL(u_1-u_0)=f(x,u_0)-ku_0-\partial^2 u_0+ku_0\ge 0,&  x\in e_j,\, j\in J\\
\sum_{j\in Inc_i} \partial_j (u_1-u_0)(v_i)\ge 0,& v_i\in V
\end{array}\right.
\]
the inequality $u_1\le u_0=u_+$ on $\G$ follows immediately by the Maximum Principle (see \cite[Prop.12]{cm}). Assuming inductively that $u_n\le u_{n-1}$, we  have  
for $x\in e_j$, $j\in J$
\begin{align*}
\cL(u_{n+1}-u_n)&=k(x)(u_{n-1}-u_n)+h(x)(e^{u_{n-1}}-e^{u_n})\\
&\ge k_1(x)e^{u_+(x)}(u_{n-1}-u_n) -k_1(x)(e^{u_{n-1}}-e^{u_n})\\
&\ge   k_1(x)(e^{u_+(x)}-e^{\xi(x)})(u_{n-1}-u_n),
\end{align*}
where $\xi(x)\in [u_n(x),u_{n-1}(x)]$. By induction,  we have $u_+\ge u_{n-1}$ and, recalling the condition at the vertices, we get
\[
\left\{\begin{array}{ll}
\cL(u_{n+1}-u_n)\ge 0&x\in e_j,\, j\in J,\\
\sum_{j\in Inc_i} \partial_j (u_{n+1}-u_{n})(v_i)\ge 0&v_i\in V;
\end{array}\right.
\]
we conclude again by the Maximum Principle that	 $u_{n+1}\le u_n$ in $\G$. We finally observe that, arguing as before, we have
\[
\left\{\begin{array}{ll}
\cL(u_-- u_{n+1})\ge k(x)(u_n-u_-)+h(x)(e^{u_{n}}-e^{u_-})\ge 0&x\in e_j,\, j\in J,\\
\sum_{j\in Inc_i} \partial_j (u_{-}-u_{n+1})(v_i)\ge 0&v_i\in V,
\end{array}\right.
\]
and therefore $u_-\le u_{n+1}$ on $\G$ for all $n$. Hence the claim \eqref{c<0;2} is proved.\\
By \cite[Prop.10]{cm} there exists a positive constant~$C$ (independent of~$n$) such that~$\|u_n\|_{H^1}\leq C$ and, in particular, $\|u_n\|_\infty\leq C$ for every~$n\in\N$.
By the first equation in~\eqref{KW} and~\eqref{c<0;2}, we deduce $\|u_n\|_{H^2}\leq C$. The Ascoli-Arzela's Theorem yields that, up to passing to a subsequence, $\{u_n\}$ converges uniformly to a function~$u\in H^1(\G)$ which is a weak solution to~\eqref{KW} with $u_-\leq u\leq u_+$.
Finally, by Remark~\ref{rmk:weak-strong}, $u$ is a classical solution to~\eqref{KW}.
\end{Proof}
In the next lemma, we show that \eqref{KW} admits a lower solution $u_-$ for any $c<0$.
\begin{Lemma}\label{c<0;L2}
If $c<0$,	there exists a lower solution $u_-$ of \eqref{KW}.
\end{Lemma}
\begin{Proof}
Set $u_-\equiv -A$ for some constant $A>0$. Then, the function~$u_-$ fulfills the Kirchhoff condition in~\eqref{KW} and also
\[\partial^2 u_-(x)-c+h(x)e^{u_-(x)}= -c+h(x)e^{-A}\ge 0\qquad  x\in e_j,\, j\in J\]
for  $A$ sufficiently large. Hence $u_-$   is a lower solution to \eqref{KW}.
\end{Proof}
\begin{Proofc}{Proof of Theorem \ref{thm:<0}}
Assume that there exists a solution $u$ of \eqref{KW}. Then, multiplying  \eqref{KW} by the test function $\phi\equiv 1$, integrating on~$\Gamma$   and taking advantage of the Kirchhoff condition and   the continuity of $u$ at the vertices, we get
\begin{align*}
&-\int_\Gamma h(x)\, dx=
\int_\Gamma(\partial u(x))^2e^{-u(x)}\, dx- c\int_\G e^{-u(x)}dx>0
\end{align*}
and therefore {\it (i)}.\par
We now assume that $\int_\Gamma h(x)\, dx<0$. Recall that, by Lemma \ref{c<0;L1} and \ref{c<0;L2},   \eqref{KW} has a solution if and only if there exists an upper solution $u_+$ to the problem. Moreover it is easy to see that, if $u_+$ is an upper solution for a given $\bar c<0$, then it is also an upper solution for any $c$ such that $ \bar c\le c<0$. Hence it follows that there exists a constant $c(h)$ with $-\infty\le c(h)\le 0$ such that \eqref{KW} admits  a solution for $c>c(h)$ and  no solution for $c<c(h)$.\\
We  show that $c(h)<0$. 
Let $m\in C^2(\G)$ be a solution of 
\begin{equation} \label{c<0;3}
\left\{\begin{array}{ll}
\partial^2 m(x) = \int_{\G} h(x)dx-h(x)&\qquad \textrm{if } x\in e_j,\, j\in J,\\
\sum_{j\in Inc_i} \partial_j m(v_i)=0&\qquad v_i\in V,
\end{array}\right.
\end{equation}
(existence of a weak solution is proved in \cite[Prop.13]{cm}, while the regularity follows by Remark \ref{rmk:weak-strong}) and   $a$ a positive constant such that 
\[ \max_{x\in\G}|e^{am(x)}-1|\le \frac{-\int_\Gamma h(x) dx}{2 \|h(x)\|_\infty}. \]
We define $b=\ln(a)$, $c=\half a\int_\G h(x)dx$ and $u_+(x)=a m(x)+b$. Then $c<0$ and  
\begin{eqnarray*}
\partial^2 u_+(x)-c +h(x)e^{u_+(x)}&=&a h(x)(e^{am(x)}-1)+ \frac{a\int_\G h(x)dx}{2}\\
&\le& a\|h(x)\|_\infty|e^{a m(x)}-1|+\frac{a\int_\G h(x)dx}{2}\leq0.
\end{eqnarray*}
Moreover, by \eqref{c<0;3}, $u_+$ is continuous and verifies the Kirchhoff condition because~$m$ enjoys the same properties. 
Hence  $u_+$ is  an upper solution and therefore we conclude that   
$$c(h)\le \frac{a}{2}\int_\G h(x)dx  <0.$$ 
We finally prove {\it (iii)}. Note that $\int h<0$ ensures $h\not\equiv 0$. \\We first show that, if $h\le 0$ in $\G$, then \eqref{KW} is solvable for any $c<0$ and therefore $c(h)=-\infty$.
Fixed $c<0$, let $m$ be a solution of \eqref{c<0;3} and choose two constants $a$, $b$ such that  $a\int_\G h(x)dx <c$ and $e^{am(x)+b}-a>0$ for $x\in\G$. We show that the function $u_+(x)=a m(x)+b$ is an upper solution of \eqref{KW}. Indeed, there holds
\begin{equation*}
\partial^2 u_+(x)-c+h(x)e^{u_+(x)}
 =a\int_\G h(x)dx -a h(x)-c + h(x)e^{a m(x)+b}\le h(x)(e^{am(x)+b}-a)\le 0
\end{equation*}
while the continuity and the Kirchhoff conditions for~$u_+$ come again from those of~$m$.
Hence $u_+$ is   an upper solution to \eqref{KW} and therefore, for any $c<0$, there exists a solution to \eqref{KW}.

Conversely, let us prove that~$c(h)=-\infty$ implies~$h\leq 0$ in~$\G$. To this end, as in~\cite[Thm2.3]{ge}, we argue by contradiction assuming that $\{h>0\}$ is not empty. For any~$c<0$, let~$u$ be a solution to~\eqref{KW} (whose existence is ensured by~$c(h)=-\infty$) and let~$\phi_c\in C^2(\G)$ be a solution to problem
\begin{equation} \label{c<0;3bis}
\left\{\begin{array}{ll}
\partial^2 \phi_c + c\phi_c =h&\quad \textrm{if } x\in e_j,\, j\in J,\\
\sum_{j\in Inc_i} \partial_j \phi_c(v_i)=0&\qquad v_i\in V
\end{array}\right.
\end{equation}
(whose existence is ensured by~\cite[Prop.10]{cm}). We claim
\begin{equation*}
\phi_c(x)\geq e^{-u(x)}>0 \qquad \forall x\in\G.
\end{equation*}
In order to prove this relation, by the Maximum Principle (\cite[Prop.12]{cm}), it suffices to prove that~$e^{-u}$ is a lower solution to~\eqref{c<0;3bis}. Actually, there holds
\[
\partial^2(e^{-u}) +c e^{-u}=e^{-u}[-\partial^2 u +|\partial u|^2+c]=e^{-u}[he^u +|\partial u|^2]\geq h;
\]
moreover, ~$e^{-u}$ is continuous and satisfies the Kirchhoff condition because~$u$ does it. Hence, our claim is proved.\\
Furthermore, we have $\lim_{c\to-\infty} c\phi_c(x)=h(x)$ for any $x\in\G$ because $(I-\partial^2)$ is a maximal monotone operator when coupled with Kirchhoff condition. Finally this property contradicts $\phi_c\geq 0$ in~$\{h>0\}$.
\end{Proofc}
\subsection{The critical case $c=c(h)$}\label{sect:ge}
\begin{Proposition}\label{critical:prop}
For $\int_\Gamma h<0$ and $c(h)>-\infty$, problem \eqref{KW} with $c=c(h)$ admits a solution.
\end{Proposition}
\begin{Proof} Note that Theorem~\ref{thm:<0}-(iii) ensures that $h$ changes sign (and obviously, $h\not\equiv 0$).
Given a decreasing sequence $\{c_k\}_{k\in \N}$ with $c(h)<c_k<0$   converging to $c(h)$ as $k\to+\infty$, we consider
\begin{equation}\label{KW_k}
\left\{\begin{array}{ll}
\partial^2 u = c_k-he^u,&\qquad \textrm{if } x\in e_j,\, j\in J,\\
\sum_{j\in Inc_i} \partial_j u(v_i)=0,&\qquad v_i\in V.
\end{array}\right.
\end{equation}
The idea is to show that a sequence of continuous solutions $u_k$ of \eqref{KW_k}, appropriately chosen, converges for $k\to \infty$ to a solution of \eqref{KW} with $c=c(h)$.
\begin{Lemma}\label{critical:L1}
For each $k\in\N$, there exist  a  lower solution $\phi_k\equiv - A\in\R$  and an upper solution $\psi_k$ to \eqref{KW_k} with $\psi_k>\phi_k$.
\end{Lemma}
\begin{Proof}
To show the existence of a lower solution, it suffices to argue as in Lemma~\ref{c<0;L2} choosing $A$ sufficiently large so that
\begin{equation}\label{crit;1}
- c_k+h(x)e^{-A}\geq - c_k - \|h\|_\infty e^{-A}=:\delta>0 . 
\end{equation}
For the upper solution,  we choose~$\psi_k$ as a solution to~\eqref{KW} with $c$ replaced by any~$\tilde c_k\in (c(h),c_k)$ (whose existence is established in Theorem~\ref{thm:<0}).

\noindent Finally, it remains to prove the inequality $\psi_k>-A$. Denoted by $\tilde x$ a minimum point of~$\psi_k$ on $\G$, we claim that $\psi_k(\tilde x)>-A$.

%

\noindent Assume first that $\tilde x\in e_j$ for some $j\in J$. The first equation in~\eqref{KW} yields:
\[
h(\tilde x)e^{\psi_k(\tilde x)}= \tilde c_k - \partial^2\psi_k(\tilde x)\leq \tilde c_k<0
\]
and, in particular,
\[h(\tilde x)<0.
\]
On the other hand,  the function~$\phi_k\equiv -A$ satisfies
\[
h(\tilde x)e^{-A}> \tilde c_k.
\]
The last three relations give: $e^{\psi_k(\tilde x)}-e^{-A}>0$, which is equivalent to~$\psi_k(\tilde x)>-A$.

\noindent Assume now $\tilde x = v_i$ for some $i\in I$ and, for later contradiction, $\psi_k(v_i)\leq -A$.
 We observe that, for any~$j\in Inc_i$, the restriction of~$\psi_{k}$ to~$e_j$ attains its minimum at~$v_i$ and, consequently, $\partial_j \psi_{k}(v_i)\geq 0$. Taking into account the Kirchhoff condition in~\eqref{KW}, we deduce
\[
\partial_j \psi_k(v_i)=0\qquad \forall j\in Inc_i.
\]
On the other hand, by \eqref{crit;1} and   the continuity of~$h$, there exists~$\eta>0$ such that
\begin{equation}\label{crit;2}
\tilde c_k+\|h\|_\infty e^{-A+\eta}<-\delta/2.
\end{equation}
Moreover, by the continuity of~$\psi_k$ and $\psi_k(v_i)\leq -A$, \eqref{crit;2} ensures
\[
\partial^2_j\psi_k(x)= \tilde c_k- h(x)e^{\psi_k(x)}\leq \tilde c_k +\|h\|_\infty e^{-A+\eta}<-\delta/2<0
\]
for any~$x\in e_j$ sufficiently near~$v_i$. In conclusion, near $v_i$, the function $\partial_j \psi_k$ is strictly decreasing with $\partial_j \psi_k(v_i)=0$ and therefore $\psi_k$ is strictly decreasing. This fact contradicts that $\psi_k$ attains its minimum at~$v_i$.
\end{Proof}
\begin{Lemma}\label{critical:L2} 
Fix $k\in\N$. The  minimum of the problem 
\begin{equation}\label{Ge1}
\inf\{{\mathcal I}_k(u):\quad u\in H^1(\Gamma),\, -A\leq u(x)\leq \psi_k(x) \quad\forall x\in\Gamma\}
\end{equation}
where
\[
{\mathcal I}_k(u):=\frac12 \int_\Gamma |\partial u|^2\, dx + c_k \int_\Gamma  u\, dx -\int_\Gamma h e^u\, dx,
\]
is attained by some function~$\bar u$  with 
\begin{equation}\label{crit;3}
-A< \bar u< \psi_k.
\end{equation}
Moreover $\bar u$ is a solution of \eqref{KW_k}.
\end{Lemma}
\begin{Proof}
Let $\{v_n\}_n$ be a minimizing sequence for ${\mathcal I}_k$. Then there holds: ${\mathcal I}_k\leq {\mathcal I}_k(-A)=c_k(-A)|\Gamma|-e^{-A}\int_\Gamma h\le C$, for some constant $C$ (independent of $k$). Moreover, we have
\begin{equation}\label{crit;4}
\begin{split}
C\geq {\mathcal I}_k(v_n)&=\frac12 \int_\Gamma |\partial v_n|^2\, dx + c_k \int_\Gamma v_n\, dx -\int_\Gamma h e^{v_n}\, dx\\
&\geq\frac12 \int_\Gamma |\partial v_n|^2\, dx + c_k \int_\Gamma  \psi_k\, dx -\|h\|_\infty\int_\Gamma e^{\psi_k}\, dx,
\end{split}
\end{equation}
where the inequality is due to the constraint $-A\leq v_n\leq \psi_k$. We deduce that $\|\partial v_n\|_2$ are uniformly bounded; on the other hand, also $\|v_n\|_\infty$ are uniformly bounded. Therefore, the sequence~$\{v_n\}_n$ is uniformly bounded in $H^1(\Gamma)$. We infer that, possibly passing to a subsequence, there exists $\bar u\in H^1(\Gamma)$ with $-A\leq \bar u\leq \psi_k$ such that: $v_n\to \bar u$ uniformly and $v_n\rightharpoonup \bar u$ weak in $H^1$. By the lower semicontinuity of ${\mathcal I}_k$, we get ${\mathcal I}_k(\bar u)\leq \liminf_n {\mathcal I}_k(v_n)$, hence $\bar u$ is minimum for \eqref{Ge1}. The
inequality \eqref{crit;3} is a consequence of the Maximum Principle.
Finally, by standard Lagrange multipliers method, we have 
\[
\frac{d}{dt}{\mathcal I}_k(\bar u+t\phi)|_{t=0}=0
\]
for any $\phi\in H^1(\G)$, from which we get \eqref{KW_weak}. Arguing as in Remark \ref{rmk:weak-strong}, we get that $\bar u$ is a strong solution to \eqref{KW_k}.
\end{Proof}
We can now conclude the proof of Proposition \ref{critical:prop}. Denote by $u_k$, $k\in \N$, a solution of \eqref{KW_k} given by Lemma \ref{critical:L2}. Assume for the moment that the sequence $\{u_k\}_k$ is bounded in $H^1(\G)$.
Hence there exists~$u\in H^1(\Gamma)$ such that, as $k\to+\infty$,  up to a subsequence, $u_k\rightharpoonup u$ in the weak topology of~$H^1(\Gamma)$ and  $u_k\to u$ uniformly. Passing to the limit in the weak formulation of \eqref{KW_k}, we get that $u$ is a weak, and therefore also a strong, solution to \eqref{KW} with $c=c(h)$.

It remains to prove that $\{u_k\}_k$ is bounded in $H^1$. To this end, fix $0<\delta<\max_\G h$, an interval $D$ inside some edge~$e_j$ such that $D\subset\{h(x)\geq \delta\}$ and a point~$\bar x\in D$; by the same arguments of~\cite[pag.743]{cl} (note that we can use \cite[Lemma2.1]{cl} because any solution of the equation in~$D$ is also a solution in a $2$-dimensional domain), we get that the $u_k$'s are uniformly bounded in $D$. Therefore, the functions $w_k(x):=u_k(x)-u_k(\bar x)$ satisfy $w_k(\bar x)=0$ and there exists $C_1>0$ such that $|u_k(\bar x)|\leq C_1$ for any $k$. Arguing as in Lemma \ref{lemma:poinc}-(i), we get: $\|w_k\|_\infty\leq |\G|^{1/2} \|\partial w_k\|_2= |\G|^{1/2} \|\partial u_k\|_2$ and, we deduce
\begin{equation}\label{crit;123}
\|u_k\|_\infty\leq |u_k(\bar x)| +\|w_k\|_\infty\leq C_1+|\G|^{1/2} \|\partial u_k\|_2.
\end{equation}
On the other hand, choosing $\phi\equiv 1$ as test function in the weak formulation of \eqref{KW_k}, we get
\begin{equation}\label{crit;5}
\int_\G he^{ u_k}dx=c_k|\G|.
\end{equation}
Since $c_k$ are negative, relations \eqref{crit;4} with $v_n=u_k$ and \eqref{crit;5} entail
\begin{eqnarray*}
C&\geq &\frac12 \int_\G |\partial u_k|^2\, dx + c_k \int_\G u_k\, dx -\int_\Gamma h e^{u_k}\, dx\geq  \frac{\|\partial u_k\|_2^2}{2} + c_k \int_\G |u_k|\, dx- c_k|\G|\\
&\geq & \frac{\|\partial u_k\|_2^2}{2}+c_kC_1|\G|+c_k|\G|^{3/2}\|\partial u_k\|_2- c_k|\G|
\end{eqnarray*}
where the last inequality is due to~\eqref{crit;123}. Hence, $\partial u_k$ are uniformly bounded in $L^2$; by~\eqref{crit;123}, the $u_k$'s are uniformly bounded in $L^\infty$ and consequently also in $H^1$.
\end{Proof}
\noindent{\bf Acknowledgements.} The second author is member of Indam-Gnampa and he has been partially supported by the fondazione Cariparo project ``Nonlinear partial differential equations: asymptotic problems and Mean Field Games''.


\end{document}